\documentclass[12pt,amstex]{article}
\usepackage{graphicx}

\usepackage{amsfonts}
\usepackage{float}
\usepackage{flafter}

\labelwidth\leftmargini\advance\labelwidth-\labelsep

\def\R{\relax\ifmmode I\!\!R\else$I\!\!R$\fi}

\def\Z{\relax\ifmmode Z\!\!\!Z\else$Z\!\!\!Z$\fi}

\def\C{\relax\ifmmode C\!\!\!\!I\else$C\!\!\!\!I$\fi}

\def\K{\relax\ifmmode I\!\!K\else$I\!\!K$\fi}

\def\N{\relax\ifmmode I\!\!N\else$I\!\!N$\fi}

\newcounter{defcounter}[section]
{\vspace{0.1cm}\begin{sloppypar}\noindent\stepcounter{defcounter}{\bfseries Definition
      \thesection.\thedefcounter}}%
{\end{sloppypar}\vspace{0.1cm}}

\newtheorem{theorem}{Theorem}[section]

\newtheorem{proposition}{Proposition}[section]

\newcommand{\proof}{{\noindent\bf Proof. }}

\newcommand{\qed}{\hfill $\square$}

\begin{document}
\thispagestyle{empty}
\begin{center}
{\Large {\bf Li-Yorke pairs of full Hausdorff dimension for some
chaotic dynamical systems}}
\end{center}
\begin{center}
J. Neunh\"auserer\\~\\
TU Clausthal ~\\
Department of Mathematics~\\
Erzstra\ss e 1~~\\
D-38678 Clausthal-Zellerfeld
Germany \\
neunchen@aol.com
\end{center}
\begin{abstract}
We show that for some simple classical chaotic dynamical
systems the set of Li-Yorke pairs has full Hausdorff dimension on
invariant sets.
\end{abstract}
\section{Introduction}
The term "chaos" in mathematical theory of dynamical systems was
used for the first time in the influential paper of Li and Yorke
\cite{[LY]}. The approach of Li and Yorke is based on the existence
of Li-Yorke pairs. These are pairs of points in the phase space that
approaches each other for some sequence of moments in the time evolution
and that remain separated for other sequences of moments. The
characteristic property of dynamical systems that are chaotic in the
sense of Li-Yorke is sensitivity to the initial condition. States that
are physically indistinguishable result in physically distinguishable
states for such systems. ~
\\Althogh there is an enormous literature on Li-Yorke chaos (see for instance \cite{[BGKM]} and \cite{[BHS]} and the
bibliography therein) there seem to be no results from a
dimensional-theoretical point of view. In this paper we ask the natural question
if for invariant sets like repellers, attractors or hyperbolic sets
of chaotic dynamical systems, Li-Yorke pairs have full Hausdorff
dimension, see section two for appropriate definitions. Our main
result in Section 5 is that this is in fact the case, if the
invariant set is self-similar or a product of self-similar sets and
the dynamics of the system is homomorphic conjugated to a full
shift, see Theorem 5.1 and 5.2 below. To prove this result we have a
look at Li-Yorke pairs in symbolic dynamics in Section 3 and
study Li-Yorke pairs in the context of iterated function systems in
Section 4. Our general result can be applied to simple
classical models of "chaotic" dynamics like the tent map, the Bakers
transformation, Smales horseshoe, and solenoid-like systems, see
Section 6. To prove that Li-Yorke pairs have full dimension for
more general hyperbolic systems could be a task for further
research.
\section{Basic Notation}
Consider a dynamical system $(X,T)$ on a complete separable metric
space $X$. A pair of points $(x,y)\in X^{2}$ is called a Li-Yorke
pair for $T$ if
\[ \liminf_{n\longmapsto\infty} d(T^{n}x,T^{n}y)=0\quad\mbox{ and }\quad \limsup_{n\longmapsto\infty}
d(T^{n}x,T^{n}y)>0,\] see \cite{[LY]}. Now given an
invariant set $\Lambda\subseteq X$, i.e. $f(\Lambda)=\Lambda$, we
define the set of Li-Yorke pairs in $\Lambda$ for $T$ by
\[LY_{T}(\Lambda)=\{(x,y)\in\Lambda^{2}|(x,y)\mbox{ is a Li-Yorke
pair}\}.\] We say that Li-Yorke pairs in $\Lambda$ have full Hausdorff
dimension for $T$ if the Hausdorff dimension of $LY_{T}(\Lambda)$
coincides with the Hausdorff dimension of $\Lambda^{2}$, i.e.
\[ \dim_{H}(LY_{T}(\Lambda))=\dim_{H}(\Lambda^{2}).\]
Recall that the Hausdorff dimension of $A\subseteq X$ is given by
\[\dim_{H}(A)=\inf\{s|H^{s}(A)=0\}=\sup\{s|H^{s}(A)=\infty\}\]
where $H^{s}(A)$ denotes the $s$-dimensional Hausdorff measure,
i.e.
\[H^{s}(A)=\lim_{\epsilon\longmapsto 0}\inf\{\sum|U_{i}|^{s}|A\subseteq \bigcup U_{i},\quad|U_{i}|<\epsilon\},\]
and $|U_{i}|$ denotes the diameter of a covering element $U_{i}$. Moreover we note that the Hausdorff dimension of a Borel probability measure $\mu$ on $X$ is defined by
\[ \dim_{H}\mu=\inf\{\dim_{H}A|\mu(A)=1\}.\]
Beside Hausdorff dimension we will use the Minkowski dimension in
our proofs for some technical reasons. Let $N_{\epsilon}(A)$ be the
smallest number of balls needed to cover $A$. The Minkowski
dimension of $A$ is given by
\[ \dim_{M}A=\lim_{\epsilon\longmapsto 0}\frac{\log
N_{\epsilon}(A)}{-\log\epsilon}, \] if the limit exists. We refer to the
book of Falconer \cite{[FA]} for an introduction to the modern
dimension theory. Given a dynamical system, invariant sets like
repellers, attractors or hyperbolic sets often have a fractal
geometry with non integer Hausdorff dimension. In this way the dimension
theory comes into the study of dynamical systems, see the book of
Pesin \cite{[PE]}. If Li-Yorke pairs have full dimension in an
invariant set for a dynamical system, this means that the measure of
chaos on this set is maximal from dimensional-theoretic point of
view.
\section{Li-Yorke pairs in symbolic dynamics}
Consider the spaces of one and two sided sequences
$\Sigma=\{1,2,\dots m\}^{\mathbb{N}}$  $\tilde\Sigma=\{1,2,\dots
m\}^{\mathbb{Z}}$ with the metric
\[ dist(\mathfrak{s},\mathfrak{t})=\sum m^{-|k|}|s_{k}-t_{k}|,\]
where $\mathfrak{s}=(s_{k})$ and $\mathfrak{t}=(t_{k})$. These are
perfect, totally disconnected and compact metric spaces. The shift
map $\sigma$ on this spaces is given by
\[\sigma((s_{k}))=(s_{k+1}).\]
For an introduction to symbolic dynamics consider for instance
\cite{[KH]}.~\\~\\ Two sequences $\mathfrak{s}$ and $\mathfrak{t}$
form a Li-Yorke pair in $\Sigma$ (or $\tilde \Sigma$) for $\sigma$
with respect to the metric $dist$, if they coincide on a sequence of
blocks with increasing length and do not coincide on one
subsequence. We construct here Li-Yorke pairs in the following way: Fix
$\mathfrak{s}\in\Sigma$ and an arbitrary sequence
$\mathfrak{N}=(N_{n})$ of natural numbers. Let the first digit of
$\mathfrak{t}\in\Sigma$ be $s_{1}$ and let the second digit by
$t_{2}=s_{2}+1$ modulo $m$. Then we choose $N_{1}$ arbitrary digits.
Next choose two digits of $\mathfrak{s}$ and one digit of
$\mathfrak{s}+1$ modulo $m$. Now we again choose $N_{2}$ arbitrary
digits and three digits from $\mathfrak{s}$ and one of
$\mathfrak{s}+1$ modulo $m$ and so on. Thus we consider subsets of
$\Sigma$ given by
\[\Sigma_{\mathfrak{N}}(\mathfrak{s})=\{\mathfrak{t}\in\Sigma|t_{k}=s_{k}\mbox{ for }k\in \{u_{i},\dots,u_{i}+i\}
,\]\[ \mbox{ and }t_{u_{i}+i+1}=s_{u_{i}+i+1}+1~\mbox{mod}~m~~\mbox{for}~i=0,\dots,\infty\}\]
where $u_{0}=1$ and $u_{i}$ is given by the recursion
$u_{i+1}=u_{i}+N_{i}+i+1$.
\begin{proposition}
A pair $(\mathfrak{s},\mathfrak{t})\in\Sigma^{2}$ with
$\mathfrak{t}\in\Sigma_{\mathfrak{N}}(\mathfrak{s})$ is a Li-Yorke
pair for $\sigma$.~\\
\end{proposition}
\proof Under the assumptions we have
\[
\lim_{i\longmapsto\infty}
d(\sigma^{u_{i}}(\mathfrak{s}),\sigma^{u_{i}}(\mathfrak{t}))\le\lim_{i\longmapsto\infty}1/m^{i}=0,\]
\[ \limsup_{i\longmapsto\infty}
d(\sigma^{u_{i}+i+1}(\mathfrak{s}),\sigma^{u_{i}+i+1}(\mathfrak{t}))\ge
1/m
\]giving the required asymptotic.\qed~\\
\\By this proposition we obviously have
\[ \Pi_{\mathfrak{N}}:=\{(\mathfrak{s},\mathfrak{t})|\mathfrak{s}\in\Sigma \mbox{ and }\mathfrak{t}\in\Sigma_{\mathfrak{N}}(\mathfrak{s})\}\subseteq
LY_{\sigma}(\Sigma),\] where $LY$ is the set of Li-Yorke pairs defined
in the preceding section. Moreover, for the two-side sequence we have
\[ \tilde\Pi_{\mathfrak{N}}:=\{(\tilde\mathfrak{s},\tilde\mathfrak{t})|\tilde \mathfrak{s}\in\Sigma \mbox{ and }
\mathfrak{t}\in\Sigma_{\mathfrak{N}}(\mathfrak{s})\}\subseteq
LY_{\sigma}(\tilde\Sigma),\] where $\mathfrak{s}$ is the part of
$\tilde\mathfrak{s}$ with positive indices. In the next section the
symbolic sets defined here will be used. In addition we will use the
natural bijection $\Sigma$ onto
$\Sigma_{\mathfrak{N}}(\mathfrak{s})$ which we denote by
$pr_{\mathfrak{s},\mathfrak{N}}$. This bijection just fills
arbitrary digits of sequences in
$\Sigma_{\mathfrak{N}}(\mathfrak{s})$ successively with a given sequence of
digits from $\Sigma$; compare with the construction of $\Sigma_{\mathfrak{N}}(\mathfrak{s})$ above. In our study of the dimension of Li-Yorke pairs
of dynamical systems that are conjugated to shift systems, the
symbolic approach will be useful in Section 5.
\section{Li-Yorke pairs for iterated function systems}
Consider a system of contracting similitudes,
$S_{i}:\mathbb{R}^{w}\longmapsto \mathbb{R}^{w}$
\[ |S_{i}x-S_{i}y|=c_{i}|x-y|\]
with $c_{i}\in(0,1)$ for $i=1,\dots m$. It is well known \cite{[HU]}
that there is a compact self-similar set $\Lambda$ with
\[ \Lambda=\bigcup_{i=1}^{m}S_{i}(\Lambda). \]
The set may be described using the projection
$\pi:\Sigma\longmapsto\Lambda$ given by
\[ \pi((s_{k}))=\lim_{n\longmapsto\infty}S_{s_{n}}\circ\dots
S_{s_{1}}(K)\] where $K$ is a compact set with $S_{i}(K)\subseteq
K$. In our study of Li-Yorke pairs for dynamical systems on $X$
we are interested in the subset of $\Lambda$ given by
\[
\Lambda_{\mathfrak{N}}(\mathfrak{s})=\pi(\Sigma_{\mathfrak{N}}(\mathfrak{s})),\]
where the set of symbols
$\Sigma_{\mathfrak{N}}(\mathfrak{s})$ was introduced in the preceding section. We use an
extension of the classical argument to proove the following result on
the dimension of $\Lambda_{\mathfrak{N}}(\mathfrak{s})$:
\begin{proposition}
Let $S_{i}:K\longmapsto K$ for $i=1\dots m$ be contracting
similitudes for some compact set $K\subseteq \mathbb{R}^{w}$ with
$S_{i}(K)\cap S_{j}(K)=\emptyset$ for $i\not=j$.\\ If
$\mathfrak{N}=(N_{n})$ is a sequence of natural numbers with
\[ \lim_{M\longmapsto\infty} M^2/\sum_{n=1}^{M}N_{n}=0,\]
 then
for all $\mathfrak{s}\in\Sigma$
\[\dim_{H}\Lambda_{\mathfrak{N}}(\mathfrak{s})=\dim_{H}\Lambda=D,\]
where $D$ is the solution of
\[ \sum_{i=1}^{m}c_{i}^{D}=1\]
for the contraction constants $c_{i}$ of the similitudes.
\end{proposition}
\proof Fix $\mathfrak{s}=(s_{k})$ and $\mathfrak{N}=(N_{n})$ throughout
the proof. Write $\bar\Lambda$ for
$\Lambda_{\mathfrak{N}}(\mathfrak{s})$, $\bar\Sigma$ for
$\Sigma_{\mathfrak{N}}(\mathfrak{s})$ and $pr$ for the bijection
from
$\Sigma$ onto $\bar\Sigma$ defined at the end of Section 3.~\\
It is well known that
$\dim_{H}\Lambda=D$, hence $\dim_{H}\bar \Lambda\le D$, see
\cite{[MO]}. For the opposite inequality we construct a Borel
probability measure $\mu$ of dimension $D$ on $\bar \Lambda$. To
this end consider the probability vector
$(c_{1}^{D},\dots,c_{m}^{D})$ on $\{1,\dots,n\}$ and the corresponding Bernoulli
measure $\nu$ on $\Sigma$, which is the infinite product of this measure. Now we map this measure onto $\bar\Sigma$
using $pr$ and further onto $\bar\Lambda$ using $\pi$, i.e.\[
\mu=\pi(pr(\nu))=\nu\circ pr^{-1}\circ \pi^{-1}.
\]
The local mass distribution princip states that
\[ \liminf_{\rho\longmapsto 0}\frac{\log\mu( B_{\rho}(x))}{\log\rho}\ge D\] for all $x\in\bar\Lambda$ implies $\dim_{H}\mu\ge D$ and hence $\dim_{H}\bar\Lambda \ge
D$. This is Proposition 2.1 of \cite{[YO]}. Hence if we prove this estimate on local dimension the proof is complete.
~\\~\\
By bijectivity of the coding map for all points $x\in \bar\Lambda$
there is a unique sequence $\mathfrak{a}=(a_{k})$ such that
$\pi(pr(\mathfrak{a}))=x$, where $pr$ is defined at the end of Section 3. We have
\[ \{x\}=\bigcap_{k=1}^{\infty} \pi(pr[a_{1},\dots,a_{k}]),\]
where $pr$ acts pointwise on cylinder sets
\[ [a_{1},\dots,a_{k}]:=\{\mathfrak{t}=(t_{i})\in\Sigma|t_{i}=a_{i}\mbox{ for } i=1,\dots k\}\]
in $\Sigma$. For further use note that $pr[a_{1},\dots,a_{k}]$ is itself a cylinder set in $\bar\Sigma$, the length of this cylinder set is $k$ plus the digits coming from the fixed sequence $\mathfrak{s}$.
\\
Given an arbitrary real $\rho>0$ we choose $k=k(\rho)$ such that
\[ d\cdot c(pr[a_{1},\dots,a_{k}])\le\rho<d\cdot
c(pr[a_{1},\dots,a_{k-1}]),\] where
$c([a_{1},\dots,a_{p}]):=c_{a_{1}}\cdot\dots\cdot c_{a_{p}}$ for all
cylinder sets and $d$ is the minimal distance of two sets in the
construction of $\bar\Lambda$, i.e. \[ d=\min_{i\not=
j}d(S_{i}(K),S_{k}(K)).\] Given another finite sequence $(\bar
a_{1},\dots,\bar a_{k})\not=(a_{1},\dots, a_{k} )$ the
contraction property of the maps $S_{i}$ implies
\[dist(\pi(pr[a_{1},\dots,a_{k}]),\pi(pr[\bar a_{1},\dots,\bar
a_{k}]))\ge d\cdot c(pr[a_{1},\dots,a_{k}])>\rho ,\] here again $[
a_{1},\dots, a_{k}]$ is a cylinder set in $\Sigma$ and
$pr[a_{1},\dots,a_{k}]$ is the corresponding cylinder set in
$\bar\Sigma$. Hence we have
\[ \bar\Lambda\cap B_{\rho}(x)\subseteq
\pi(pr[a_{1},\dots,a_{k}])\] and
\[ \mu(B_{\rho}(x))\le \mu(\pi(pr[a_{1},\dots,a_{k}]))=(c([a_{1},\dots,a_{k}]))^{D}\]
\[
=\frac{(c([a_{1},\dots,a_{k}]))^{D}}{(c([pr(a_{1},\dots,a_{k}]))^{D}}(c(pr[a_{1},\dots,a_{k}]))^{D}\le\frac{(c([a_{1},\dots,a_{k}]))^{D}}{(c(pr[a_{1},\dots,a_{k}]))^{D}}
d^{-D}\rho^{D}\] by the construction of the measure $\mu$ and the
choice of $k$. Now taking logarithm this yields
\[\log\mu(B_{\rho}(x))\le D(\log\rho-\log
d-\log\frac{c(pr[a_{1},\dots,a_{k}])}{c([a_{1},\dots,a_{k}])}) \le
D(\log\rho-\log d-\sharp(k)\log\underline{c}), \] where
$\underline{c}=\min_{i=1}^{m} c_{i}$ and $\sharp(k)$ is the length
of the cylinder set $pr[a_{1},\dots,a_{k}]$ minus $k$, the length of the
cylinder set $[a_{1},\dots,a_{k}]$. Now dividing by $\log\rho$ and using
the definition of $k$ we obtain,
\[\frac{\log\mu(B_{\rho}(x))}{\log\rho}\ge
D+D(-\frac{\log d}{\log\rho}-\frac{\sharp(k)\log\underline{c}}{\log
d+\log(c(pr[a_{1},\dots,a_{k-1}]))})\]\
 \[\ge D+D(-\frac{\log
d}{\log\rho}-\frac{\sharp(k)\log\underline{c}}{\log
d+((k-1)+\sharp(k-1))\log\overline{c} })
,\]
where $\overline{c}=\max_{i=1}^{m} c_{i}$ . We have
$\lim_{\rho\longmapsto 0}k(\rho)=\infty$, hence it remains to show
that
\[\lim_{k\longmapsto\infty}\frac{\sharp(k)}{k}=0.\]
Given $k$ choose $M(k)$ such that
\[ \sum_{n=1}^{M(k)-1}N_{n}<k\le \sum_{n=1}^{M(k)}N_{n}.\]
By the definition of $\bar \Sigma$ and the map $pr$ we have
\[\sharp(k)\le \sum_{v=1}^{M(k)}(v+1)<M(k)^{2},\]
hence
\[\frac{\sharp(k)}{k}\le \frac{M(k)^{2}}{\sum_{n=1}^{M(k)-1}N_{n}}.\]
By the assumption on $(N_{n})$ the righthand side goes to zero with
$k\longmapsto \infty$. This completes the proof. \qed ~\\~\\Using
general results in dimensions theory we may go one step further and
show that the projection of the set $\Pi_{\mathfrak{N}}$ has full
dimension in $\Lambda^{2}$. This is the image of Li-Yorke pairs for
the shift map in the symbolic space, compare with Section 3.
We will use the following fact:
\begin{theorem}
Let $F$ be a subset of $\mathbb{R}^{n}$ and $E$ be a subset of $\mathbb{R}^{k}$ with $k<n$. Let $L_{x}$ be the $n-k$ dimensional affine linear subspace of $\mathbb{R}^{n}$ given by the translation $x\in E$. If $\dim_{H}(F\cap L_{x})\ge t$ for all $x\in E$ than $\dim_{H}F\ge t+\dim_{H}E$.
\end{theorem}
This is Corollary 7.1 of Falconer \cite{[FA]} based on the work of Marstrand \cite{[MA]}. With the help of this theorem we get:
\begin{proposition}
Under the assumptions of Proposition 4.1 we have
\[ \dim_{H}S=\dim_{H}\Lambda^{2}=2\dim_{H}\Lambda\] where
$S=\pi(\Pi_{\mathfrak{N}})=\{(x,y)|x\in\Lambda~,~
y\in\Lambda_{\mathfrak{N}}(\pi^{-1}(x))\}$.
\end{proposition}
\proof Since Hausdorff and Minkowski dimension of $\Lambda$ coincide
we have $\dim_{H}\Lambda^{2}=2\dim_{H}\Lambda$, see Corollary 7.4 of
\cite{[FA]}. Obviously $S\subseteq \Lambda^{2}$, hence $\dim_{H}S\le
2\dim_{H}\Lambda$. On the other hand Proposition 4.1 implies
$\dim_{H}(S\cap
\Lambda_{\mathfrak{N}}(\pi^{-1}(x)))=\dim_{H}\Lambda$ for all
$x\in\Lambda$. By Theorem 4.1 we get $\dim_{H}S\ge2\dim_{H}\Lambda$, which completes the proof. \qed ~\\~\\We remark at the end of this
section that the results proved here for $\mathbb{R}^{w}$ remain
true on complete separable metric spaces of finite multiplicity
which have the Besicovitsh property, compare Appendix I of
\cite{[PE]}. The techniques we have used apply in the general setting.
\section{Li-Yorke pairs of full dimension for systems conjugated to a shift}
In this section we state and prove our main results on the Hausdorff
dimension of Li-Yorke pairs of dynamical systems conjugated to a shift and
having a self-similar invariant set. The results are consequences of
Proposition 3.1, Proposition 4.1 and Proposition 4.2 below.
\begin{theorem}
Let $f:\mathbb{R}^{w}\longmapsto \mathbb{R}^{w}$ be a dynamical
system with a compact invariant set $\Lambda$. If $(\Lambda,f)$ is
homomorphic conjugated to a one-sided full shift
$(\Sigma,\sigma)$ and $\Lambda$ is self-similar then the Li-Yorke pairs
in $\Lambda$ have full Hausdorff dimension for $f$.
\end{theorem}
\proof By the assumption of the theorem there is a hom\"omorphism
$\pi:\Sigma\longmapsto \Lambda$ with
\[ f\circ\pi=\pi\circ f. \quad (\star) \]
If $(s,t)\in\Sigma^{2}$ is a Li-Yorke pair for $\sigma$ then by definition
\[ \liminf_{n\longmapsto\infty} d(\sigma^{n}s,\sigma^{n}t)=0\mbox{ and } \limsup_{n\longmapsto\infty}
d(\sigma s,\sigma t)>0.\] By continuity this implies
\[ \liminf_{n\longmapsto\infty} d(\pi(\sigma^{n}s),\pi(\sigma^{n}t))=0\mbox{ and } \limsup_{n\longmapsto\infty}
d(\pi(\sigma^{n} s),\pi(\sigma^{n} t))>0\] and using $(\star)$
\[\liminf_{n\longmapsto\infty} d(f^{n}(\pi(s)),f^{n}\pi(t))=0\mbox{ and } \limsup_{n\longmapsto\infty}
d(f^{n}(\pi(s)),f^{n}(\pi(t)))>0,\] which means
$(\pi(s),\pi(t))\in\Sigma^{2}$ is a Li-Yorke pair of $f$. Since
$\Pi_{\mathfrak{N}}\subseteq LY_{\sigma}(\Sigma)$, see Section 3,
we get
 \[ S=\pi(\Pi_{\mathfrak{N}})\subseteq LY_{f}(\Lambda).\]
Furthermore, since $\Lambda$ is self-similar and fulfils the
condition of Proposition 4.1, by the properties of the coding $\pi$
we have by Proposition 4.2 $\dim_{H}S=\dim_{H}\Lambda^{2}$ and hence
$\dim_{H} LY_{f}(\Lambda)=\dim_{H}\Lambda^{2}$ concluding the proof.
\qed ~\\~\\For dynamical systems which have a product structure we
obtain the following result:
\begin{theorem}
Let $f:\mathbb{R}^{w}\longmapsto \mathbb{R}^{w}$ be a dynamical
system with a compact invariant set $\Lambda$. If $(\Lambda,f)$ is
hom\"oormorphic conjugated to a two-side full shift
$(\tilde\Sigma,\sigma)$ and the coding map is a product with two
self-similar images, $\Lambda=\Lambda_{1}\times\Lambda_{2}$, then
the Li-Yorke pairs in $\Lambda$ have full dimension for $f$.
\end{theorem}
\proof By the same argument as in the proof of the Theorem 5.1 we have
 $S=\pi(\tilde\Pi_{\mathfrak{N}})\subseteq
LY_{f}(\Lambda)$. Using $\pi=(\pi_{1},\pi_{2})$ we have
\[ S=\{(\pi_{1}(s^{+}),\pi_{2}(s^{-}),\pi_{1}(t^{+}),\pi_{2}(
t^{-})|\tilde\mathfrak{s}=(s^{+},s^{-})\in\tilde\Sigma \quad
\tilde\mathfrak{t}=(t^{+},t^{-})\in\tilde\Sigma\mbox{ with }
t^{+}\in\Sigma_{\mathfrak{N}}(\mathfrak{s}^{+})\}.\] Since $\Lambda$
is self-similar and fulfils the condition of Proposition 4.1, by the
properties of the coding $\pi_{1}$ we have
\[ \dim_{H}\{\pi(t^{+})|t^{+}\in\Sigma_{\mathfrak{N}}(\mathfrak{s}^{+})\}=\dim_{H}\Lambda_{1}.\]
By the argument used in the proof of Proposition 4.2 this implies
\[ \dim_{H}S\ge 2\dim_{H}\Lambda_{1}+2\dim_{H}\Lambda_{2}, \]
see again Corollary 7.12 of \cite{[FA]} and \cite{[MA]}. On the
other hand, since Minkowski and Hausdorff dimensions of self-similar
sets $\Lambda_{1},\Lambda_{2}$ coincides, we have
\[\dim_{H}S\le
dim_{H}\Lambda=2\dim_{H}\Lambda_{1}+2\dim_{H}\Lambda_{2}\] by
Corollary 7.4 of \cite{[FA]}. This concludes the
proof.\qed~\\~\\
The last section of this paper is devoted to examples.
\section{Examples}
In this section we consider four classical examples of chaotic
dynamical systems, namely the tent map in dimension one, the skinny
Backers transformation and a linear horseshoe in dimension two and
linear solenoid-like systems in dimension three. All these systems
have self-similar invariant sets with dynamics conjugated to a full
shift on two symbols.~\\~\\
First consider the expansive tent map, see \cite{[KH]},
$t:\mathbb{R}\longmapsto\mathbb{R}$ given by
\[ t(x)=a-2a|x-1/2|, \]
where $a>1$. The map has an invariant repeller $\Lambda$ which is
given by the iterated function system
\[ T_{1}x=\frac{1}{2a}x, \qquad T_{2}x=1-\frac{1}{2a}x.\]
The dynamical system $(\Lambda,t)$ is conjugated to a one-sided
shift on two symbols via the coding homomorphism induced by the
iterated function system:
\[ \pi(\mathfrak{s})=\lim_{n\longmapsto\infty} T_{s_{1}}\circ
T_{s_{2}}\circ\dots\circ T_{s_{n}}([0,1]). \] By Theorem 4.1 the Li-Yorke
pairs in $\Lambda$ have full Hausdorff dimension for $t$.~\\~\\
Now consider the skinny Backers transformation, see \cite{[NE]},
$b:[0,1]^{2}\longmapsto [0,1]^{2}$ given by
\[ b(x,y)=\{\begin{array}{cc}(\beta_{1} x, 2y)~~~~~~~~~~~~~~~~~\mbox{ if } y\le 1/2\\(1-\beta_{2}+\beta_{2}x,1-2y)\mbox{ if } y> 1/2
\end{array}\]
for $\beta_{1},\beta_{2}\in (0,1)$ with $\beta_{1}+\beta_{2}<1$. The
map has an attractor given by $\Lambda\times [0,1]$, where $\Lambda$
is given by the iterated function system
\[ T_{1}x=\beta_{1}x,\quad T_{2}x=1-\beta_{2}+\beta_{2}x.\]
The system $(\Lambda\times [0,1],b)$ is homomorphic conjugated to
a two-sided full shift on two symbols via $\pi=(\pi_{1},\pi_{2})$
where $\pi_{1}$ is given by the iterated function system and
$\pi_{2}$ is just the map coming from dyadic expansion. The assumptions of
Theorem 4.2 are fulfilled and we again have Li-Yorke pairs of full
Hausdorff dimension.~\\~\\
Next consider Smale's horseshoe $h:[0,1]^{2}\longmapsto
\mathbb{R}^{2}$, see \cite{[SM]}, fulfilling
\[ h(x,y)=\{\begin{array}{cc}(\beta x, \tau y)~~~~~~~~~~~~~~~~~\mbox{ if } y\le 1/\tau  \\(-\beta x+1,-\tau
x+\tau )~~~~~~~\mbox{ if } y>1-1/\tau
\end{array}
\]
on horizontal strips. Here we assume $\beta\in(0,1/2)$ and $\tau>2$.
The map may be extended to a diffeomorphism of $\mathbb{R}^{2}$
using stretching and folding of the middle strip $1/\tau
>y>1-1/\tau$. The hyperbolic invariant set
\[ \Lambda=\bigcap_{n=-\infty}^{\infty}b([0,1]^{2}) \]
is given by $\Lambda=\Lambda_{1}\times\Lambda_{2}$ where
$\Lambda_{1}$ is given by the iterated function system
\[ T_{1}x=\beta x,\qquad T_{2}x=-\beta x+1\]
and $\Lambda_{2}$ is given by the iterated function system
\[ G_{1}y=\frac{1}{\tau} y,\qquad G_{2}y=-\frac{1}{\tau} y+1.\]
To this iterated function system there corresponds a homomorphic
coding $\pi=(\pi_{1},\pi_{2})$ of
$(\Lambda_{1}\times\Lambda_{2},b)$. Again the assumptions of Theorem
4.2 are fulfilled and we have Li-Yorke pairs of full
Hausdorff dimension.\\~\\
Our last example is a solenoid like system, see \cite{[NE2]}, given
by $s:[0,1]^{3}\longmapsto[-1,1]^{3}$ given by
$b:[0,1]^{2}\longmapsto [0,1]^{2}$ given by
\[ s(x,y,z)=\{\begin{array}{cc}(\beta_{1} x, \beta_{1}y, 2z)~~~~~~~~~~~~~~~~~~~~~~~~~~~~~\mbox{ if } z\le 1/2\\(1-\beta_{2}+\beta_{2}x,1-\beta_{2}+\beta_{2}y,1-2z)\mbox{ if } z> 1/2
\end{array}\]
for $\beta_{1},\beta_{2}\in (0,1)$ with $\beta_{1}+\beta_{2}<1$. By
exactly the same argument we used in the case of the skinny Bakers
transformation we see that the Li-Yorke pairs have full Hausdorff
dimension on the attractor $\Lambda$ for the map $s$.

\end{document}